\documentclass{elsart}

\usepackage{amsmath}
\usepackage{amscd,amssymb}
\usepackage{leftidx}

\begin{document}

\begin{frontmatter}

 \title{The $r$-Bell numbers}
 \author{Istv\'an Mez\H{o}\thanksref{a}}
 \thanks[a]{Present address: University of Debrecen, H-4010, Debrecen, P.O. Box 12, Hungary}
 \address{Department of Applied Mathematics and Probability Theory, Faculty of Informatics, University of Debrecen, Hungary}
 \ead{mezo.istvano@inf.unideb.hu}
 \ead[url]{http://www.inf.unideb.hu/valseg/dolgozok/mezoistvan/mezoistvan.html}

\begin{abstract}The notion of generalized Bell numbers has appeared in several works but there is no systematic treatise on this topic. In this paper we fill this gap. We discuss the most important combinatorial, algebraic and analytic properties of these numbers which generalize the similar properties of the Bell numbers. Most of these results seem to be new. It turns out that in a paper of Whitehead these numbers appeared in a very different context. In addition, we introduce the so-called $r$-Bell polynomials.
\end{abstract}

\begin{keyword}Bell numbers \sep $r$-Bell numbers \sep Stirling numbers \sep $r$-Stirling numbers
\MSC 11B73
\end{keyword}
\end{frontmatter}


\newcommand{\NN}{\mathbb{N}}
\newcommand{\ZZ}{\mathbb{Z}}
\newcommand{\QQ}{\mathbb{Q}}
\newcommand{\im}{\mbox{\upshape Im}}

\newcommand{\stirlings}[2]{\genfrac\{\}{0pt}{}{#1}{#2}}
\newcommand{\stirlingf}[2]{\genfrac[]{0pt}{}{#1}{#2}}

\newtheorem{Theorem}{Theorem}[section]
\newtheorem{Proposition}[Theorem]{Proposition}
\newtheorem{Lemma}[Theorem]{Lemma}
\newtheorem{Corollary}[Theorem]{Corollary}
\newtheorem{Conjecture}[Theorem]{Conjecture}
\theoremstyle{definition}
\newtheorem{Example}[Theorem]{Example}
\newtheorem{Remark}[Theorem]{Remark}
\newtheorem{Problem}[Theorem]{Problem}

\section{Introduction}

The $B_n$ Bell number \cite{Comtet} counts the partitions of a set with $n$ elements. Since $\stirlings{n}{k}$, for a fixed $n$ and $k$, enumerates the number of partitions of a set with $n$ elements consisting $k$ disjoint, nonempty sets, we get immediately that $B_n$ can be given by the sum
\begin{equation}
B_n=\sum_{k=0}^n\stirlings{n}{k}.\label{Belldef}
\end{equation}
The numbers $\stirlings{n}{k}$ are called Stirling numbers of the second kind. The $n$-th Bell polynomial is
\[B_n(x)=\sum_{k=0}^n\stirlings{n}{k}x^k.\]
These numbers and polynomials have many interesting properties and appear in several combinatorial identities. A comprehensive paper is \cite{Comtet}.

There is a more general version of Stirling numbers, the so-called $r$-Stirling numbers. The number $\stirlings{n}{k}_r$, for fixed $n\ge k\ge r$, enumerates the  partitions of a set of $n$ elements into $k$ nonempty, disjoint subsets such that the first $r$ elements are in distinct subsets. A systematic treatment on the $r$-Stirling numbers is given in \cite{Broder}, and a different approach is described in \cite{Carlitz1,Carlitz2}. According to \eqref{Belldef}, it seems to be natural to define the numbers
\begin{equation}
B_{n,r}=\sum_{k=0}^n\stirlings{n+r}{k+r}_r.\label{rBelldef}
\end{equation}
(It is obvious that $B_n=B_{n,0}$, because $\stirlings{n}{k}=\stirlings{n}{k}_0$ by the definitions.)

The very first question is on the meaning of the $r$-Bell numbers. By \eqref{rBelldef}, $B_{n,r}$ is the number of the partitions of a set with $n+r$ element such that the first $r$ elements are in distinct subsets in each partition.

The name of $r$-Stirling numbers suggests the name for the numbers $B_{n,r}$: we call them as $r$-Bell numbers, and the name of the polynomials
\[B_{n,r}(x)=\sum_{k=0}^n\stirlings{n+r}{k+r}_rx^k\]
will be $r$-Bell polynomials. Thus $B_{n,r}=B_{n,r}(1)$.

As far as we know, this paper is the first one fully devoted to the $r$-Bell numbers, although Carlitz \cite{Carlitz1,Carlitz2} defined these numbers and proved some identities for them. His original notation was
\[B_{n,r}=B(n,r).\]

\section{Example and tables}

The following example illuminates again the meaning of the $r$-Bell numbers. By definition,
\[B_{2,2}=\stirlings{4}{2}_2+\stirlings{4}{3}_2+\stirlings{4}{4}_2.\]
$\stirlings{4}{2}_2$ counts the partitions of $4$ element into $2$ subsets such that the first $2$ element are in distinct subsets:
\[\{1,3,4\},\{2\}\quad;\quad\{1\},\{2,3,4\}\quad;\quad\{1,3\},\{2,4\}\quad;\quad\{1,4\},\{2,3\}.\]
$\stirlings{4}{3}_2$ belongs to the partitions
\[\{1\},\{2\},\{3,4\}\quad;\quad\{1,3\},\{2\},\{4\}\quad;\quad\{1,4\},\{2\},\{3\}\quad;\]
\[\{1\},\{2,3\},\{4\}\quad;\quad\{1\},\{2,4\},\{3\}.\]
Finally, $\stirlings{4}{4}_2$ equals to the number of partitions of $4$ elements into $4$ subsets (and trivially, the first two elements are in distinct subsets again):
\[\{1\},\{2\},\{3\},\{4\}.\]
That is,
\[B_{2,2}=\stirlings{4}{2}_2+\stirlings{4}{3}_2+\stirlings{4}{4}_2=4+5+1=10\]
is really the number of all partitions of the set $\{1,2,3,4\}$ such that the first two element are in distinct subsets.

\vspace{1cm}
\begin{center}\textbf{The first $r$-Bell numbers}\end{center}\vspace{5mm}
\begin{tabular}{r|cccccccc}&$n=0$&$n=1$&$n=2$&$n=3$&$n=4$&$n=5$&$n=6$&\\\hline
$r=0$&1&1&2&5&15&52&203\\
$r=1$&1&2&5&15&52&203&877\\
$r=2$&1&3&10&37&151&674&3263\\
$r=3$&1&4&17&77&372&1915&10481\\
$r=4$&1&5&26&141&799&4736&29371\\
$r=5$&1&6&37&235&1540&10427&73013\\
$r=6$&1&7&50&365&2727&20878&163967
\end{tabular}

\vspace{1cm}
\begin{center}\textbf{The first $r$-Bell polynomials}\end{center}
\begin{eqnarray*}
B_{0,r}(x)&=&1\\
B_{1,r}(x)&=&x+r\\
B_{2,r}(x)&=&x^2+(2r+1)x+r^2\\
B_{3,r}(x)&=&x^3+(3r+3)x^2+(3r^2+3r+1)x+r^3\\
B_{4,r}(x)&=&x^4+(4r+6)x^3+(6r^2+12r+7)x^2+\\&&(4r^3+6r^2+4r+1)x+r^4
\end{eqnarray*}

\section{Generating functions}

We start to derive the properties of $r$-Bell numbers. First of all, the generating functions are determined.

\begin{Theorem}The exponential generating function for the $r$-Bell polynomials is
\[\sum_{n=0}^\infty B_{n,r}(x)\frac{z^n}{n!}=e^{x(e^z-1)+rz}.\]
\end{Theorem}

\textit{Proof. }Broder \cite{Broder} gave the double generating function of $r$-Stirling numbers
\[\sum_{n=0}^\infty\left(\sum_{k=0}^n\stirlings{n+r}{k+r}x^k\right)\frac{z^n}{n!}=e^{x(e^z-1)+rz}.\]
Although he did not use the inner sum in the past of his paper, they are exactly our polynomials. We note that this identity is remarked in \cite[eq. (3.19)]{Carlitz1}
\hfill\qed

\begin{Corollary}The $r$-Bell polynomials can be expressed by the usual Bell ones.
\[B_{n,r}(x)=\sum_{k=0}^nr^k\binom{n}{k}B_{n-k}(x).\]
\end{Corollary}

\textit{Proof. }In the following steps we use the well known identity for ordinary Bell polynomials:
\[B_n(x+y)=\sum_{k=0}^n\binom{n}{k}B_k(x)B_{n-k}(y).\]
Then the exponential generating function of the $r$-Bell polynomials can be transformed as
\begin{eqnarray*}
&&\sum_{n=0}^\infty B_{n,r}\left(\frac{x}{2}+\frac{x}{2}\right)\frac{z^n}{n!}=e^{\left(\frac{x}{2}+\frac{x}{2}\right)(e^z-1)+rz}\\
&=&\left(\sum_{n=0}^\infty r^n\frac{z^n}{n!}\right)\left[\left(\sum_{n=0}^\infty B_n\left(\frac{x}{2}\right)\frac{z^n}{n!}\right)\left(\sum_{n=0}^\infty B_n\left(\frac{x}{2}\right)\frac{z^n}{n!}\right)\right]\\
&=&\left(\sum_{n=0}^\infty r^n\frac{z^n}{n!}\right)\left(\sum_{n=0}^\infty B_n\left(\frac{x}{2}+\frac{x}{2}\right)\frac{z^n}{n!}\right)\\
&=&\sum_{n=0}^\infty\left(\sum_{k=0}^n r^k\binom{n}{k}B_{n-k}(x)\right)\frac{z^n}{n!},
\end{eqnarray*}
by the Cauchy product. Comparing the coefficients, the Corollary follows.\hfill\qed

We remark that the non-polynomial version was proven by Carlitz \cite[eq. (3.18)]{Carlitz1}.

In order to determine the ordinary generating function we need some other notions. The falling factorial of a given real number $x$ is denoted and defined by
\begin{equation}
x^{\underline{n}}=x(x-1)(x-2)\cdots(x-n+1),\quad(n=1,2,\dots)\label{fallfact}
\end{equation}
and $(x)^{\underline{0}}=1$, while the rising factorial (a.k.a. Pochhammer symbol) is
\begin{equation}
(x)_n\equiv x^{\overline{n}}=x(x+1)(x+2)\cdots(x+n-1)\quad(n=1,2,\dots)\label{poch}
\end{equation}
with $(x)_0=1$. It is obvious that $(1)_n=n!$. Fitting our notations to the theory of hypergeometric functions introduced immediately, we apply the notation $(x)_n$ instead of $x^{\overline{n}}$. The next transformation formula holds
\begin{equation}
x^{\underline{n}}=(-1)^n(-x)_n.\label{pochtrans}
\end{equation}
The hypergeometric function (or hypergeometric series) is defined by the series
\[\leftidx{_p}{F}{_q}\left(\left.\begin{tabular}{llll}$a_1,$&$a_2,$&$\dots,$&$a_p$\\$b_1,$&$b_2$,&$\dots,$&$b_q$\end{tabular}\right|t\right)=\sum_{k=0}^\infty\frac{(a_1)_k(a_2)_k\cdots(a_p)_k}{(b_1)_k(b_2)_k\cdots(b_q)_k}\frac{t^k}{k!}.\]
The ordinary generating function can be given with this function.

\begin{Theorem}The $r$-Bell polynomials have the generating function
\[\sum_{n=0}^\infty B_{n,r}(x)z^n=\frac{-1}{rz-1}\,\frac{1}{e^x}\,\leftidx{_1}{F}{_1}\left(\left.\begin{tabular}{llll}$\;\,\frac{rz-1}{z}$\\$\frac{rz+z-1}{z}$\end{tabular}\right|x\right).\]
\end{Theorem}

\textit{Proof. }It is known \cite{Broder} that for the Stirling numbers
\[\sum_{n=0}^\infty\stirlings{n}{m}_rz^n=\frac{z^m}{(1-rz)(1-(r+1)z)\cdots(1-mz)}\quad(m\ge r\ge 0).\]
This can be rewritten as
\[\sum_{n=m}^\infty\stirlings{n+r}{m+r}_rz^n=\frac{z^m}{(1-rz)(1-(r+1)z)\cdots(1-(m+r)z)}.\]
We transform the denominator using the falling factorial:
\[(1-rz)(1-(r+1)z)\cdots(1-(m+r)z)\]
\[=\frac{(1-z)(1-2z)\cdots(1-(m+r)z)}{(1-z)(1-2z)\cdots(1-(r-1)z)}=\frac{z^{m+1}\left(\frac{1}{z}\right)^{\underline{m+r+1}}}{\left(\frac{1}{z}\right)^{\underline{r}}}.\]
Hence
\[\sum_{k=m}^\infty\stirlings{k+r}{m+r}_rz^k=\frac{1}{z}\left(\frac{1}{z}\right)^{\underline{r}}\frac{1}{\left(\frac{1}{z}\right)^{\underline{m+r+1}}}.\]
Rule \eqref{pochtrans} and the definitions \eqref{fallfact}-\eqref{poch} give that
\[\left(\frac{1}{z}\right)^{\underline{m+r+1}}=(-1)^{m+r+1}\left(-\frac{1}{z}\right)_{m+r+1}\]
\[=(-1)^{m+r+1}\left(-\frac{1}{z}\right)_{r+1}\left(-\frac{1}{z}+r+1\right)_m.\]
Consequently,
\[\sum_{n=m}^\infty\stirlings{n+r}{m+r}_rz^n=\frac{1}{z}\frac{\left(\frac{1}{z}\right)^{\underline{r}}}{\left(-\frac{1}{z}\right)_{r+1}}\frac{(-1)^{m+r+1}}{\left(\frac{rz+z-1}{z}\right)_m}.\]
Since
\[\frac{\left(\frac{1}{z}\right)^{\underline{r}}}{\left(-\frac{1}{z}\right)_{r+1}}=(-1)^r\frac{z}{rz-1},\]
we get that
\[\sum_{n=m}^\infty\stirlings{n+r}{m+r}_rz^n=\frac{-1}{rz-1}\frac{(-1)^m}{\left(\frac{rz+z-1}{z}\right)_m}.\]
Let us multiplication with $x^m$ and take a summation for $m$ over the non-negative integers,
\[\sum_{n=0}^\infty B_{n,r}(x)z^n=\frac{-1}{rz-1}\sum_{m=0}^\infty\frac{(-x)^m}{\left(\frac{rz+z-1}{z}\right)_m}=\frac{-1}{rz-1}\,\leftidx{_1}{F}{_1}\left(\left.\begin{tabular}{llll}$\;\;\;\;1$\\$\frac{rz+z-1}{z}$\end{tabular}\right|-x\right).\]
Finally, we apply Kummer's formula \cite[p. 505.]{Abramowitz}
\[e^{-x}\,\leftidx{_1}{F}{_1}\left(\left.\begin{tabular}{llll}$a$\\$b$\end{tabular}\right|x\right)=\,\leftidx{_1}{F}{_1}\left(\left.\begin{tabular}{llll}$b-a$\\$\;\;\;b$\end{tabular}\right|-x\right)\]
with $b=\frac{rz+z-1}{z}$ and $a=\frac{rz-1}{z}$.
\hfill\qed

\section{Basic recurrences}

In an earlier paper of the author \cite{Mezo}, the polynomials $B_{n,r}(x)$ were introduced because of a very different reason. These functions were used to study the properties of the $r$-Stirling numbers and some necessary properties of them were proven in that paper. We repeat those results without proof.

\begin{Theorem}We have the next recursive identities:
\begin{eqnarray*}
B_{n,r}(x)&=&x\left(\frac{d}{dx}B_{n-1,r}(x)+B_{n-1,r}(x)\right)+rB_{n-1,r}(x),\\
e^xx^rB_{n,r}(x)&=&x\frac{d}{dx}\left(e^xx^rB_{n-1,r}(x)\right).
\end{eqnarray*}
Moreover, all roots of $B_{n,r}(x)$ are real and negative.
\end{Theorem}

A straightforward corollary that for a fixed $r$ the constant term of the $n$-th polynomial is $r^n$:
\[B_{n,r}(0)=r^n\]
and that the derivative of an $r$-Bell polynomial is determined by the relation
\[\frac{\partial}{\partial x}B_{n,r}(x)=\frac{B_{n+1,r}(x)}{x}-\frac{rB_{n,r}(x)}{x}-B_{n,r}(x).\]

The identity
\[\stirlings{n+r}{k+r}_r=\stirlings{n+r}{k+r}_{r-1}-(r-1)\stirlings{n-1+r}{k+r}_{r-1}\]
was proven in \cite{Broder} and implies the recursive relation
\[B_{n,r}(x)=B_{n,r-1}(x)-(r-1)B_{n-1,r-1}(x).\]

Finally, we cite Carlitz's identities \cite[eq. (3.22-3.23)]{Carlitz1}:
\begin{eqnarray}
B_{n+m,r}&=&\sum_{j=0}^m\stirlings{m+r}{j+r}_rB_{n,r+j},\label{Carlitzeq}\\
B_{n,r+m}&=&\sum_{j=0}^m(-1)^{m-j}\stirlingf{m+r}{j+r}_rB_{n+j,r}.\nonumber
\end{eqnarray}
Here $\stirlingf{n}{m}_r$ is an $r$-Stirling number of the first kind (see \cite{Broder,Carlitz1,Carlitz2}).

\section{Dobinski's formula}

The Bell numbers are involved in Dobinski's nice formula \cite{Chowla,Dobinski,GKP,Pitman}:
\[B_n=\frac{1}{e}\sum_{k=0}^\infty\frac{k^n}{k!}.\]
Our goal is to generalize this identity to our case.

\begin{Theorem}[Dobinski's formula]The $r$-Bell polynomials satisfy the i\-den\-ti\-ty
\[B_{n,r}(x)=\frac{1}{e^x}\sum_{k=0}^\infty\frac{(k+r)^n}{k!}x^k.\]
Consequently, the $r$-Bell numbers are generated by
\[B_{n,r}=\frac{1}{e}\sum_{k=0}^\infty\frac{(k+r)^n}{k!}.\]
\end{Theorem}

\textit{Proof. }The $r$-Stirling numbers for a fixed $n$ (and $r$) have the 'horizontal' generating function \cite{Broder}
\[(x+r)^n=\sum_{k=0}^n\stirlings{n+r}{k+r}_rx^{\underline{k}},\]
whence, for an arbitrary integer $m$,
\[\frac{(m+r)^n}{m!}=\sum_{k=0}^m\stirlings{n+r}{k+r}_r\frac{1}{(m-k)!}.\]
In the next step we multiply both sides with $x^m$ and sum from $m=0$ to $\infty$. Then
\[\sum_{m=0}^\infty\frac{(m+r)^n}{m!}x^m=\sum_{m=0}^\infty\sum_{k=0}^m\stirlings{n+r}{k+r}_r\frac{x^m}{(m-k)!}=e^x\left(\sum_{k=0}^n\stirlings{n+r}{k+r}_rx^k\right).\]
\hfill\qed

We can determine some interesting sum with the aid of $r$-Bell numbers. For example, we know from the second paragraph that $B_{2,2}=10$, so
\[\frac{1}{e}\sum_{k=0}^\infty\frac{(k+2)^2}{k!}=10.\]

\section{An integral representation}

In 1885, Ces\`aro \cite{Cesaro} found an amazing integral representation of the Bell numbers (see also \cite{Becker,Callan}):
\[B_n=\frac{2n!}{\pi e}\,\im\int_0^\pi e^{e^{e^{i\theta}}}\sin(n\theta)d\theta.\]
It is not hard to deduce the '$r$-Bell version'.

\begin{Theorem}The $r$-Bell numbers have the integral representation as
\[B_{n,r}=\frac{2n!}{\pi e}\,\im\int_0^\pi e^{e^{e^{i\theta}}}e^{re^{i\theta}}\sin(n\theta)d\theta.\]
\end{Theorem}

\textit{Proof. }In \cite{Carlitz1} we find that
\begin{equation}
k!\stirlings{n+r}{k+r}_r=\sum_{j=0}^k(-1)^{k-j}\binom{k}{j}(j+r)^n.\label{stirlingeq}
\end{equation}
In \cite{Callan} the next equality appears:
\begin{equation}
\im\int_0^\pi e^{je^{i\theta}}\sin(n\theta)d\theta=\frac{\pi}{2}\frac{j^n}{n!}.\label{sinint}
\end{equation}
Unifying the equations \eqref{stirlingeq}-\eqref{sinint}, we get that
\begin{eqnarray*}
&&\frac{\pi}{2}\frac{1}{n!}\stirlings{n+r}{k+r}_r=\frac{1}{k!}\sum_{j=0}^k(-1)^{k-j}\binom{k}{j}\im\int_0^\pi e^{(j+r)e^{i\theta}}\sin(n\theta)d\theta\\
&=&\frac{1}{k!}\im\int_0^\pi\left[\sum_{j=0}^k(-1)^{k-j}\binom{k}{j} \left({e^{e^{i\theta}}}\right)^j\right]e^{re^{i\theta}}\sin(n\theta)d\theta\\
&=&\im\int_0^\pi\frac{\left(e^{e^{i\theta}}-1\right)^k}{k!}e^{re^{i\theta}}\sin(n\theta)d\theta,
\end{eqnarray*}
whence
\[\sum_{k=0}^\infty\stirlings{n+r}{k+r}_r=\frac{2n!}{\pi}\im\int_0^\pi\left(\sum_{k=0}^\infty\frac{\left(e^{e^{i\theta}}-1\right)^k}{k!}\right)e^{re^{i\theta}}\sin(n\theta)d\theta,\]
and the result follows.\hfill\qed

\begin{Remark}The imaginary part of the above integral can be calculated with a bit of effort. The result is
\[B_{n,r}=\frac{2n!}{\pi e}\int_0^\pi e^{e^{cos\theta}\cos\sin\theta+r\cos\theta}\cdot\]
\[\cdot\left[\cos(e^{\cos\theta}\sin\sin\theta)\sin(r\sin\theta)+\sin(e^{\cos\theta}\sin\sin\theta)\cos(r\sin\theta)\right]\sin(n\theta)d\theta.\]
Without the $r$-Bell numbers in background, the evaluation of this integral seems to be impossible\dots

Citing the general version of Dobinski's formula we find the compelling identity
\[\sum_{k=0}^\infty\frac{(k+r)^n}{k!}=\frac{2n!}{\pi}\,\im\int_0^\pi e^{e^{e^{i\theta}}}e^{re^{i\theta}}\sin(n\theta)d\theta.\]
\end{Remark}

\section{Hankel transformation and log-convexity}

Since
\[e^t\sum_{n=0}^\infty B_{n,r}(x)\frac{t^n}{n!}=e^{x(e^t-1)+(r+1)t},\]
Cauchy's product immediately implies the next

\begin{Theorem}\label{thm:binomtransform}The $r$-Bell polynomials satisfy the relations
\begin{eqnarray*}
B_{n,r+1}(x)&=&\sum_{k=0}^n\binom{n}{k}B_{k,r}(x),\\
B_{n,r}(x)&=&\sum_{k=0}^n\binom{n}{k}(-1)^{n-k}B_{k,r+1}(x).
\end{eqnarray*}
\end{Theorem}

An interesting corollary is connected with the notion of Hankel transform. The $H$ Hankel matrix \cite{Layman} of an integer sequence $(a_n)$ is
\[H=\begin{pmatrix}a_0&a_1&a_2&a_3&\cdots\\a_1&a_2&a_3&a_4&\cdots\\a_2&a_3&a_4&a_5&\cdots\\\vdots&\vdots&\vdots&\vdots&\ddots\end{pmatrix},\]
while the Hankel matrix of order $n$, denoted by $h_n$, is the upper-left submatrix of $H$ of size $n\times n$. The Hankel transform of the sequence $(a_n)$ is again a sequence formed by the determinants of the matrices $h_n$.

A notable result of Aigner and Lenard \cite{Aigner,Lenard} is that the Hankel transform of the sequence of Bell numbers is $(1!,1!2!,1!2!3!,\dots)$, that is, for any fixed $n$,
\[\begin{vmatrix}B_0&B_1&B_2&\cdots&B_n\\B_1&B_2&B_3&\cdots&B_{n+1}\\\vdots&\vdots&\vdots&&\vdots\\B_n&B_{n+1}&B_{n+2}&\cdots&B_{2n}\end{vmatrix}=\prod_{i=0}^n i!\]

We can determine the Hankel transform of $r$-Bell numbers easily. To reach this aim, we recall an other notion. If $(a_n)$ is a sequence, then its binomial transform $(b_n)$ is defined by the relation
\[b_n=\sum_{k=0}^n\binom{n}{k}(-1)^{n-k}a_k,\]
while the inverse transform is
\[a_n=\sum_{k=0}^n\binom{n}{k}b_k.\]
See the paper \cite{Riordan} on these transformations, for instance. A useful theorem of Layman \cite{Layman} states that any integer sequence has the same Hankel transform as its binomial transform and vice versa. Then Theorem \ref{thm:binomtransform} yield the

\begin{Corollary}The $r$-Bell numbers have the Hankel transform
\[\begin{vmatrix}B_{0,r}&B_{1,r}&B_{2,r}&\cdots&B_{n,r}\\B_{1,r}&B_{2,r}&B_{3,r}&\cdots&B_{n+1,r}\\\vdots&\vdots&\vdots&&\vdots\\B_{n,r}&B_{n+1,r}&B_{n+2,r}&\cdots&B_{2n,r}\end{vmatrix}=\prod_{i=0}^n i!\]
\end{Corollary}

Now we give an other consequence of Theorem \ref{thm:binomtransform}. It is known \cite{LiuWang} that the Bell numbers are log-convex, that is,
\[B_{n-1}B_{n+1}\ge B_n^2\quad(n\ge1).\]
A proposition (also known as Davenport-P\'olya theorem) in \cite{LiuWang} states that the binomial transform preserves the log-convexity. This implies that
\[B_{n-1,r}B_{n+1,r}\ge B_{n,r}^2\quad(n\ge1)\]
holds for all $r>0$, too.

Professor J. Cigler \cite{Cigler} calculated more general identities with respect to Hankel determinants involving not only $r$-Bell numbers but polynomials. We cite his unpublished results here.

Let $d(n,k)=\det(B_{i+j+k,r}(x))_{i,j=0}^{n-1}$. This is a bit modification comparing to our general Hankel matrix $h_n$. Cigler's results are the followings:
\[d(n,0)=x^{\binom{n}{2}}\prod_{k=0}^{n-1}k!,\]
and
\[d(n,1)=x^{\binom{n}{2}}\prod_{k=0}^{n-1}k!\sum_{k=0}^n\binom{n}{k}x^k(r)_{n-k}.\]

\section{Some occurrences of the $r$-Bell numbers}

Surprisingly, the $r$-Bell numbers were turned up in a table of Whitehead's paper \cite{Whitehead}. In his table, the $(n,i)$-entry is denoted by $b_{n,i}$ and it is the sum of the coefficients of the polynomial $x^i(x)_{n-i}$ with respect to the so-called complete graph base. A more detailed description on this graph theoretical notion can be found in the paper \cite{Whitehead} and the references therein.

It is easy to see that our $r$-Bell numbers are exactly the entries of that table, more exactly,
\begin{equation}
B_{n,r}=b_{n+r,n}\quad(n\ge1).\label{whitrec}
\end{equation}
From this observation we get straightaway the next identity.

\begin{Theorem}We have for all $n\ge 1$ that
\[B_{n+1,r}=rB_{n,r}+B_{n,r+1}.\]
\end{Theorem}

\textit{Proof. }According to \cite{Whitehead}, the entries $b_{n,i}$ satisfies the recurrence
\[(n-i)b_{n,i}+b_{n+1,i}=b_{n+1,i+1}.\]
Equality \eqref{whitrec} implies the statement. On the other hand, this theorem is a special case of \eqref{Carlitzeq} but it is worth to give a different viewpoint.
\hfill\qed

We note that the 'row sum' in the table of Whitehead can be expressed by the $r$-Bell numbers, too.
\[\sum_{i=1}^nb_{n,i}=\sum_{i=1}^nB_{i,n-i}.\]

The observation \eqref{whitrec} gives that the $r$-Bell numbers have meaning in the theory of chromatic polynomials.

An other occurrence is the following.
The $r$-Bell numbers come from a problem on the maximum of $r$-Stirling numbers (see \cite{Mezo}). The author proved there that all roots of the polynomial $B_{n,r}(x)$ are non-positive. This implies that
\[\stirlings{n}{k}_r^2\ge\stirlings{n}{k+1}_r\stirlings{n}{k-1}_r,\]
which is an important relation -- for example -- in the theory of combinatorial sequences. In addition, the maximizing index of $r$-Stirling numbers of the second kind can be expressed approximately by the $r$-Bell numbers \cite{Mezo}. Namely,
\[\left|K-\left(\frac{B_{n+1,r}}{B_{n,r}}-(r+1)\right)\right|<1,\]
where $K$ is the parameter, for which
\[\stirlings{n+r}{K}_r\ge\stirlings{n+r}{k}_r\]
for all $k=r,r+1,\dots,n+r$.

We remark, that (beside these and \cite{Carlitz1,Carlitz2}), there is an other paper in which the $r$-Bell numbers appear. C. B. Corcino \cite{Corcino} deals with the asymptotic properties of these numbers.

\begin{center}
 \textbf{Acknowledgement}
\end{center}
I thank Professor Cigler for his suggestions and results on Hankel deteminants of $r$-Bell polynomials.

\end{document}